\magnification=1200


\catcode`\À=\active \defÀ{\`A}    \catcode`\à=\active \defà{\`a} 
\catcode`\Â=\active \defÂ{\^A}    \catcode`\â=\active \defâ{\^a} 
\catcode`\Æ=\active \defÆ{\AE}    \catcode`\æ=\active \defæ{\ae}
\catcode`\Ç=\active \defÇ{\c C}   \catcode`\ç=\active \defç{\c c}
\catcode`\È=\active \defÈ{\`E}    \catcode`\è=\active \defè{\`e} 
\catcode`\É=\active \defÉ{\'E}    \catcode`\é=\active \defé{\'e} 
\catcode`\Ê=\active \defÊ{\^E}    \catcode`\ê=\active \defê{\^e} 
\catcode`\Ë=\active \defË{\"E}    \catcode`\ë=\active \defë{\"e} 
\catcode`\Î=\active \defÎ{\^I}    \catcode`\î=\active \defî{\^\i}
\catcode`\Ï=\active \defÏ{\"I}    \catcode`\ï=\active \defï{\"\i}
\catcode`\Ô=\active \defÔ{\^O}    \catcode`\ô=\active \defô{\^o} 
\catcode`\Ù=\active \defÙ{\`U}    \catcode`\ù=\active \defù{\`u} 
\catcode`\Û=\active \defÛ{\^U}    \catcode`\û=\active \defû{\^u} 
\catcode`\Ü=\active \defÜ{\"U}    \catcode`\ü=\active \defü{\"u} 

\catcode`\ =\active \def { }

\hsize=11.25cm    
\vsize=18cm       
\parindent=12pt   \parskip=5pt     

\hoffset=.5cm   
\voffset=.8cm   

\pretolerance=500 \tolerance=1000  \brokenpenalty=5000

\catcode`\@=11

\font\eightrm=cmr8         \font\eighti=cmmi8
\font\eightsy=cmsy8        \font\eightbf=cmbx8
\font\eighttt=cmtt8        \font\eightit=cmti8
\font\eightsl=cmsl8        \font\sixrm=cmr6
\font\sixi=cmmi6           \font\sixsy=cmsy6
\font\sixbf=cmbx6

\font\tengoth=eufm10 
\font\eightgoth=eufm8  
\font\sevengoth=eufm7      
\font\sixgoth=eufm6        \font\fivegoth=eufm5

\skewchar\eighti='177 \skewchar\sixi='177
\skewchar\eightsy='60 \skewchar\sixsy='60

\newfam\gothfam           \newfam\bboardfam

\def\tenpoint{
  \textfont0=\tenrm \scriptfont0=\sevenrm \scriptscriptfont0=\fiverm
  \def\rm{\fam\z@\tenrm}
  \textfont1=\teni  \scriptfont1=\seveni  \scriptscriptfont1=\fivei
  \def\oldstyle{\fam\@ne\teni}\let\old=\oldstyle
  \textfont2=\tensy \scriptfont2=\sevensy \scriptscriptfont2=\fivesy
  \textfont\gothfam=\tengoth \scriptfont\gothfam=\sevengoth
  \scriptscriptfont\gothfam=\fivegoth
  \def\goth{\fam\gothfam\tengoth}
  
  \textfont\itfam=\tenit
  \def\it{\fam\itfam\tenit}
  \textfont\slfam=\tensl
  \def\sl{\fam\slfam\tensl}
  \textfont\bffam=\tenbf \scriptfont\bffam=\sevenbf
  \scriptscriptfont\bffam=\fivebf
  \def\bf{\fam\bffam\tenbf}
  \textfont\ttfam=\tentt
  \def\tt{\fam\ttfam\tentt}
  \abovedisplayskip=12pt plus 3pt minus 9pt
  \belowdisplayskip=\abovedisplayskip
  \abovedisplayshortskip=0pt plus 3pt
  \belowdisplayshortskip=4pt plus 3pt 
  \smallskipamount=3pt plus 1pt minus 1pt
  \medskipamount=6pt plus 2pt minus 2pt
  \bigskipamount=12pt plus 4pt minus 4pt
  \normalbaselineskip=12pt
  \setbox\strutbox=\hbox{\vrule height8.5pt depth3.5pt width0pt}
  \let\bigf@nt=\tenrm       \let\smallf@nt=\sevenrm
  \normalbaselines\rm}

\def\eightpoint{
  \textfont0=\eightrm \scriptfont0=\sixrm \scriptscriptfont0=\fiverm
  \def\rm{\fam\z@\eightrm}
  \textfont1=\eighti  \scriptfont1=\sixi  \scriptscriptfont1=\fivei
  \def\oldstyle{\fam\@ne\eighti}\let\old=\oldstyle
  \textfont2=\eightsy \scriptfont2=\sixsy \scriptscriptfont2=\fivesy
  \textfont\gothfam=\eightgoth \scriptfont\gothfam=\sixgoth
  \scriptscriptfont\gothfam=\fivegoth
  \def\goth{\fam\gothfam\eightgoth}
  
  \textfont\itfam=\eightit
  \def\it{\fam\itfam\eightit}
  \textfont\slfam=\eightsl
  \def\sl{\fam\slfam\eightsl}
  \textfont\bffam=\eightbf \scriptfont\bffam=\sixbf
  \scriptscriptfont\bffam=\fivebf
  \def\bf{\fam\bffam\eightbf}
  \textfont\ttfam=\eighttt
  \def\tt{\fam\ttfam\eighttt}
  \abovedisplayskip=9pt plus 3pt minus 9pt
  \belowdisplayskip=\abovedisplayskip
  \abovedisplayshortskip=0pt plus 3pt
  \belowdisplayshortskip=3pt plus 3pt 
  \smallskipamount=2pt plus 1pt minus 1pt
  \medskipamount=4pt plus 2pt minus 1pt
  \bigskipamount=9pt plus 3pt minus 3pt
  \normalbaselineskip=9pt
  \setbox\strutbox=\hbox{\vrule height7pt depth2pt width0pt}
  \let\bigf@nt=\eightrm     \let\smallf@nt=\sixrm
  \normalbaselines\rm}

\tenpoint

\def\pc#1{\bigf@nt#1\smallf@nt}         \def\pd#1 {{\pc#1} }

\catcode`\;=\active
\def;{\relax\ifhmode\ifdim\lastskip>\z@\unskip\fi
\kern\fontdimen2  -1.2 \fontdimen3 \string;}

\catcode`\:=\active
\def:{\relax\ifhmode\ifdim\lastskip>\z@\unskip\fi\penalty\@M\ \fi\string:}

\catcode`\!=\active
\def!{\relax\ifhmode\ifdim\lastskip>\z@
\unskip\fi\kern\fontdimen2  -1.1 \fontdimen3 \string!}

\catcode`\?=\active
\def?{\relax\ifhmode\ifdim\lastskip>\z@
\unskip\fi\kern\fontdimen2  -1.1 \fontdimen3 \string?}

\catcode`\«=\active 
\def«{\raise.4ex\hbox{%
 $\scriptscriptstyle\langle\!\langle$}}

\catcode`\»=\active 
\def»{\raise.4ex\hbox{%
 $\scriptscriptstyle\rangle\!\rangle$}}

\frenchspacing

\def\raggedbottom{\topskip 10pt plus 36pt\r@ggedbottomtrue}

\def\pointir{\unskip . --- \ignorespaces}

\def\Medbreak{\vskip-\lastskip\medbreak}

\long\def\th#1 #2\enonce#3\endth{
   \Medbreak\noindent
   {\pc#1} {#2\unskip}\pointir{\it #3}\smallskip}

\def\demonstration{\vskip-\lastskip\smallskip\noindent
 {\it Démonstration} : }

\def\decale#1{\smallbreak\hskip 28pt\llap{#1}\kern 5pt}
\def\decaledecale#1{\smallbreak\hskip 34pt\llap{#1}\kern 5pt}
\def\puce{\smallbreak\hskip 6pt{$\scriptstyle\bullet$}\kern 5pt}

\def\eqalign#1{\null\,\vcenter{\openup\jot\m@th\ialign{
\strut\hfil$\displaystyle{##}$&$\displaystyle{{}##}$\hfil
&&\quad\strut\hfil$\displaystyle{##}$&$\displaystyle{{}##}$\hfil
\crcr#1\crcr}}\,}

\catcode`\@=12

\showboxbreadth=-1  \showboxdepth=-1

\mathcode`A="7041 \mathcode`B="7042 \mathcode`C="7043 \mathcode`D="7044
\mathcode`E="7045 \mathcode`F="7046 \mathcode`G="7047 \mathcode`H="7048
\mathcode`I="7049 \mathcode`J="704A \mathcode`K="704B \mathcode`L="704C
\mathcode`M="704D \mathcode`N="704E \mathcode`O="704F \mathcode`P="7050
\mathcode`Q="7051 \mathcode`R="7052 \mathcode`S="7053 \mathcode`T="7054
\mathcode`U="7055 \mathcode`V="7056 \mathcode`W="7057 \mathcode`X="7058
\mathcode`Y="7059 \mathcode`Z="705A

\def\qp{{\bf Q}_p}
\def\azero#1{A_0(#1)}
\def\azeroo#1{A_0(#1)_0}

\def\Kbar{{\overline K}}

\def\Xprim{{X'}}

\def\XKbar{{X{}_{\Kbar}}}

\def\Atilde{{\widetilde A}}
\def\Atildeetoile{{\widetilde A}^\times}
\def\Abaretoile{{\overline A}^\times}
\def\Ctilde{{\widetilde C}}
\def\Utilde{{\widetilde U}}

\def\Xtilde{{\widetilde X}}
\def\xtilde{{\tilde x}}
\def\ctilde{{\tilde c}}
\def\vtilde{{\tilde v}}
\def\psitilde{{\tilde\psi}}
\def\XKtilde{X_\Ktilde}
\def\Xktilde{X_\ktilde}
\def\toroxtilde{\tor^{\Orond_{\tilde X}}}
\def\Letoile{L^\times}

\def\xibar{\overline{\xi}}
\def\deltabar{\overline{\delta}}

\def\cbar{\overline{c}}

\def\Ketoile{K^{\times}}
\def\Ktilde{{\widetilde K}}
\def\ktilde{{\tilde k}}
\def\Ktildeetoile{{\Ktilde{}^{\times}}}
\def\Kbaretoile{{\Kbar{}^{\times}}}
\def\Kprim{{K'}}

\def\kprim{{k'}}
\def\det{\mathop{\hbox{d{\'e}t}}\nolimits}
\def\Gal{\mathop{\rm Gal}\nolimits}
\def\Pic{\mathop{\rm Pic}\nolimits}
\def\Br{\mathop{\rm Br}\nolimits}
\def\Hom{\mathop{\rm Hom}\nolimits}
\def\Ext{\mathop{\rm Ext}\nolimits}
\def\Card{\mathop{\rm Card}\nolimits}
\def\Spec{\mathop{\rm Spec}\nolimits}
\def\Supp{\mathop{\rm Supp}\nolimits}

\def\P{\mathord{\bf P}}
\def\Z{\mathord{\bf Z}}
\def\Ogoth{{\goth O}}
\def\ogoth{{\goth o}}
\def\ogothetoile{{{\goth o}^\times}}
\def\ogothprim{{{\goth o}'}}
\def\ogothtilde{{\tilde{\goth o}}}
\def\ogothtildemult{\ogothtilde^\times}

\def\ocycle{\hbox{$0$-cycle}}

\def\hfl#1#2#3{\smash{\mathop{\hbox to#3{\rightarrowfill}}\limits
^{\textstyle#1}_{\textstyle#2}}}
\def\gfl#1#2#3{\smash{\mathop{\hbox to#3{\leftarrowfill}}\limits
^{\textstyle#1}_{\textstyle#2}}}

\def\pafl#1#2#3{\llap{$\textstyle #1$}
\left\Vert\vbox to#3{}\right.\rlap{$\textstyle #2$}}

\def\qed{\raise -2pt\hbox{\vrule\vbox to 10pt{\hrule width 4pt
                 \vfill\hrule}\vrule}}
\def\cqfd{\unskip\penalty 500\quad\qed\medbreak}

\def\phi{\varphi}

\def\numero{n$^{\rm o}$}

\def\Ker{\mathop{\rm Ker}\nolimits}
\def\Im{\mathop{\rm Im}\nolimits}
\def\droite#1{\,\hfl{#1}{}{8mm}\,}
\def\versbas#1{\vfl{#1}{}{5mm}}

\def\diagram#1{\def\normalbaselines{\baselineskip=0pt\lineskip=5pt}
\matrix{#1}}
\def\vfl#1#2#3{\llap{$\textstyle #1$}
\left\downarrow\vbox to#3{}\right.\rlap{$\textstyle #2$}}

\newcount\refno 

\long\def\ref#1:#2<#3>{                                        
\global\advance\refno by1\par\noindent                              
\llap{[{\bf\number\refno}]\ }{#1} \pointir{\it #2} #3\goodbreak }

\def\citer#1(#2){[{\bf\number#1}\if#2\empty\relax\else,\ #2\fi]}

\def\fleche{\rightarrow}

\newcount\numerodesection
\def\section#1{\bigbreak
 {\bf\number\numerodesection.\ \ #1}\nobreak\medskip
 \advance\numerodesection by1}

\newcount\numeroderemarque
\def\remarque{\advance\numeroderemarque by1\smallbreak
{\it Remarque\/}\ \number\numeroderemarque~:}

\newcount\formuleno
\def\numeroter{\global\advance\formuleno by1
 \leqno{(\number\formuleno)}}
\def\formule(#1){{\rm (\number#1)}}

\def\div{\mathop{\rm div}\nolimits}
\def\cl{\mathop{\rm cl}\nolimits}
\def\rg{\mathop{\rm rg}\nolimits}
\def\pr{\mathop{\rm pr}\nolimits}
\def\somme{\mathop{\smash{\raise 2pt\hbox{$\sum$}}}\limits}

\def\Erond{{\cal E}}
\def\Frond{{\cal F}}
\def\Irond{{\cal I}}
\def\Lrond{{\cal L}}
\def\Nrond{{\cal N}}
\def\Orond{{\cal O}}
\def\Trond{{\cal T}}
\def\Zrond{{\cal Z}}

\def\tor{{\Trond}\!\!or}
\def\zero{\{0\}}
\def\long{\mathop{\rm long}\nolimits}
\def\ev{\mathop{\hbox{\rm év}}\nolimits}

\newbox\bibbox
\setbox\bibbox\vbox{\bigbreak
\centerline{{\pc RÉFÉRENCES BIBLIOGRAPHIQUES}}
 
\ref {\pc BARLOW} (R.):
A simply connected surface of general type with $p_g=0$, 
<Inventiones math.\ {\bf 79} (1985), 293--301.>
\newcount\barlowsimply  \global\barlowsimply=\refno

\ref {\pc BARLOW} (R.):
Rational equivalence of zero-cycles for some more surfaces with
$p_g=0$, 
<Inventiones math.\ {\bf 79} (1985), 303--308.>
\newcount\barlowrational  \global\barlowrational=\refno

\ref {\pc BERTHELOT} (P.), {\pc GROTHENDIECK} (A.) et {\pc ILLUSIE} (L.):
Théorie des intersections et théorème de
Riemann-Roch, 
<Séminaire de Géométrie Algébrique du
Bois-Marie, 1966--67, Springer, Lecture Notes {\bf 225}, 1971.>
\newcount\sgavi \global\sgavi=\refno

 
\ref {\pc BLOCH} (S.):
On the Chow groups of certain rational surfaces,
<Ann.\ sci.\ École norm.\ sup.\ (4) {\bf 14} (1981), no. 1, 41--59.>
\newcount\blochchow  \global\blochchow=\refno

\ref {\pc BLOCH} (S.):
Lettre à Colliot-Thélène,
<12 mai~1983.>
\newcount\lettre \global\lettre=\refno

\ref {\pc BOREL} (A.) et {\pc SERRE} (J-P.):
Le théorème de Riemann-Roch (d'après Grothendieck),
<Bulletin de la Soc.\ math.\ de France {\bf 86} (1958), 97--136.>
\newcount\borelserre \global\borelserre=\refno

\ref {\pc BOSCH} (S.), {\pc LÜTKEBOHMERT} (W.) et {\pc RAYNAUD} (M.):
Néron models, 
<Springer, 1990.>
\newcount\neron  \global\neron=\refno

\ref {\pc BOURBAKI} (N.):
Algèbre homologique,
<Masson, 1980.>
\newcount\bourbaki \global\bourbaki=\refno

\ref {\pc COLLIOT}-{\pc THÉLÈNE} (J-L.):
Hilbert's theorem $90$ for $K_2$, with applications to Chow groups of
rational surfaces, 
<Inventiones math.\ {\bf 71} (1983), 1--20.>
\newcount\hilbertxc \global\hilbertxc=\refno


\ref {\pc COLLIOT}-{\pc THÉLÈNE} (J-L.) et {\pc SANSUC} (J-J.):
On the Chow groups of certain rational
surfaces : a sequel to a paper of S.~Bloch, 
<Duke math.\ jour.\ {\bf 48} (1981), 421--427.>
\newcount\ctssequel \global\ctssequel=\refno

\ref {\pc COLLIOT}-{\pc THÉLÈNE} (J-L.) et {\pc SANSUC} (J-J.):
La descente sur les variétés rationnelles II,
<Duke math.\ jour.\ {\bf 54} (1987), 375--492.>
\newcount\ctsdescente \global\ctsdescente=\refno

\ref {\pc COOMBES} (K. R.) et {\pc MUDER} (D. J.):
Zero-cycles on del Pezzo surfaces over local fields, 
<Journal of Algebra {\bf 97} (1985), 438--460.>
\newcount\coomuzero \global\coomuzero=\refno

\ref{\pc DALAWAT} (C. S.):
Groupe des classes de $0$-cycles sur les surfaces
rationnelles définies sur un corps local,
<Thèse, Université de Paris-Sud, Orsay, 1993.>
\newcount\these  \global\these=\refno

\ref{\pc DALAWAT} (C. S.):
Le groupe de Chow d'une surface de Châtelet sur un corps local,
<Indag. mathem. N.S. {\bf 11} (2) (2000), 173--185, {\tt math.AG/0302156}.>
\newcount\chatelet  \global\chatelet=\refno

\ref{\pc DALAWAT} (C. S.):
The Chow group of a del Pezzo surface over a local field,
<{\tt math.AG/0302260}.>
\newcount\delpezzo  \global\delpezzo=\refno

\ref {\pc FULTON} (W.):
Intersection Theory,
<Springer, 1984.>
\newcount\fulton  \global\fulton=\refno

\ref {\pc GROTHENDIECK} (A.) et {\pc DIEUDONNÉ} (J.):
Éléments de Géométrie algébrique,
<Publications mathématiques de l'IHÉS.>
\newcount\ega \global\ega=\refno

\ref {\pc GROTHENDIECK} (A.) et {\pc DEMAZURE} (M.):
Structure des schémas en groupes réductifs,
<Séminaire de Géométrie Algébrique du Bois-Marie, 1962--64,
Springer, Lecture Notes {\bf 153}, 1970.>
\newcount\sgaiii \global\sgaiii=\refno

\ref{\pc MANIN} (Yu.\ I.):
Lectures on the $K$-functor in
algebraic geometry, 
<Russian Mathematical Surveys {\bf 24} (5) (1969), 1--89.>
\newcount\manin  \global\manin=\refno

\ref {\pc SERRE} (J-P.):
Cohomologie galoisienne,
<Springer, Lecture Notes {\bf 5}, cinquième édition, 1994.>
\newcount\serrecg  \global\serrecg=\refno

} 

\centerline{\bf Le groupe de Chow d'une surface rationnelle sur un corps
 local}

\vskip4mm

\centerline{Chandan Singh {\pc DALAWAT}}

\vskip2cm

Soient $p$ un nombre premier et $K$ une extension finie du corps
$\qp$.  Désignons par $\Kbar$ une clôture algébrique de~$K$, par
$\Ktilde$ l'extension maximale non ramifiée de $K$ dans $\Kbar$, par
$\ogoth$ (resp.\ $\ogothtilde$) l'anneau des entiers de $K$
(resp.~$\Ktilde$) et par $k$ (resp.~$\ktilde$) le corps résiduel de
$K$ (resp.~$\Ktilde$).

Soit $X$ un $\ogoth$-schéma régulier, projectif et plat de dimension
relative~2 dont la fibre générique $X_K=X\times_\ogoth K$ est {\it
potentiellement birationnelle} au plan projectif\/ $\P_2$ : on demande
que la $\Kbar$-surface $\XKbar=X_K\times_K \Kbar$ soit birationnelle à
$\P_{2,\Kbar}$; d'habitude, on dit alors que la surface $X_K$ est
rationnelle.

Notons $\azero{X_K}$ le groupe de Chow de $0$-cycles sur~$X_K$
\citer\fulton() et $\azeroo{X_K}$ le sous-groupe noyau de l'homomorphisme degré
$\deg:\azero{X_K}\fleche\Z$.  Ces groupes ont été étudiés par Bloch
\citer\blochchow(), par  Colliot-Thélène et Sansuc \citer\ctssequel()
et par Colliot-Thélène \citer\hilbertxc().  On sait ainsi que le
groupe $A_0(X_K)_0$ est fini et qu'il s'annule lorsque $X$ est lisse
sur~$\ogoth$ \citer\hilbertxc(). 

Dans cette Note, {\it on se propose de décrire un procédé pour
calculer le groupe $\azero{X_K}$ lorsque l'application naturelle\/
$\Pic X_\Ktilde\fleche\Pic X_\Kbar$ (où
$\XKtilde=X_K\times_K \Ktilde$) est un isomorphisme}.  Suivant une
idée émise par Bloch \citer\lettre() et reprise dans
\citer\these(), ce but est atteint au~\numero~5.

Décrivons ce procédé dans le cas où toutes les composantes
irréductibles rendues réduites de la fibre fermée
$X_k=X\times_\ogoth k$ de~$X$ sont lisses sur~$k$.  Soit $S$
l'ensemble (fini) des composantes irréductibles (rendues réduites) de
$X_{\ktilde}$ ; prolongeons l'action du groupe $\Gal(\ktilde|k)$
sur~$S$ au $\Z$-module libre $\Z^S$ par linéarité.  Notons
$\Hom_k(\Z^S,\Z)$ le groupe des homomorphismes
$\Gal(\ktilde|k)$-équivariants.

Pour une courbe $C$ tracée sur une composante irréductible $Y$ de
$X_k$ et pour une composante irréductible (rendue réduite) $Z\in S$ de
$X_{\ktilde}$, considérons le nombre
$$
\deg_{\Ctilde}(\Orond_{\Xtilde}(Z)|_{\Ctilde})\in\Z,
\quad\hbox{où\ }\;\Xtilde=X\times_\ogoth\ogothtilde\;\hbox{\ et\ } 
\;\Ctilde=C\times_k \ktilde. 
\numeroter
\newcount\lenombre \global\lenombre=\formuleno
$$ 
Ces nombres induisent une application $\Pic Y\fleche\Hom_k(\Z^S,\Z)$,
à savoir celle qui envoie la classe de $C$ dans $\Pic Y$ sur
l'homomorphisme $\Z^S\fleche\Z$ qui associe à $e_Z$ le
nombre~\formule(\lenombre), où $(e_Z)_{Z\in S}$ est la base canonique
de $\Z^S$.  En faisant varier $Y$, on obtient une application
$$
\bigoplus_Y \Pic Y \fleche \Hom_k(\Z^S,\Z)
\quad(Y:\ \hbox{composantes irréductibles de\ } X_k).
\numeroter
\newcount\defbx \global\defbx=\formuleno
$$
Il est clair que le conoyau de \formule(\defbx) est aisément
calculable dans des cas concrets: voir l'{\it Exemple} placé à la fin. 

Nous allons construire un isomorphisme de $\azero{X_K}$ avec le
conoyau de~\formule(\defbx).  Ce résultat fournit donc un moyen
pratique pour calculer~$\azero{X_K}$.
\smallskip
{\it Cette Note est tirée d'une thèse\/}~\citer\these() {\it préparée
sous la direction de Jean-Louis Colliot-Thélène ; je le remercie
vivement pour m'avoir mis entre les mains\/}~\citer\lettre(), {\it ainsi
que pour ses conseils et pour ses encouragements.  Je tiens à
remercier le Ministère français des Affaires étrangères pour une
bourse doctorale et l'Institut des Hautes Études scientifiques de
Bures-sur-Yvette pour deux séjours pendant lesquels la présente
rédaction a été réalisée.}

\section{Notations et hypothèses}

Les notations suivantes auront cours dans tout l'article : 

{\openup1\jot 
\halign{---\thinspace\thinspace#&$#$&$_{#}$:&\quad#&
\qquad---\thinspace\thinspace#&$#$&$_{#}$:&\quad#\hfil\cr
&p&&\multispan5\quad un nombre premier,\hfil\cr
&K&&\multispan5\quad une extension finie de $\qp$,\hfil\cr
&\Kbar&&\multispan5\quad une clôture algébrique de~$K$,\hfil\cr
&\Ktilde&&\multispan5\quad l'extension maximale non ramifiée de $K$
dans $\Kbar$,\hfil\cr 
&\ogoth&&\multispan5\quad l'anneau des entiers de $K$,\hfil\cr
&\ogothtilde&&\multispan5\quad  l'anneau des entiers de $\Ktilde$,\hfil\cr
&k&&\multispan5\quad  le corps résiduel de $K$,\hfil\cr
&\ktilde&&\multispan5\quad le corps résiduel de~$\Ktilde$,\hfil\cr
&X&&\multispan5\quad un $\ogoth$-schéma régulier, projectif et plat de
dimension relative~2,\hfil\cr
&X&K& $=X\times_\ogoth K$,&&v&&la valuation normalisée $\Ketoile\fleche\Z$,\cr
&X&\Kbar& $=X_K\times_K \Kbar$,&&\pi&&une uniformisante de $K$,\cr
&X&\Ktilde& $=X_K\times_K \Ktilde$,&&i&&l'immersion fermée $X_k\fleche X$,\cr
&X&k& $=X\times_\ogoth k$,&&j&&l'immersion ouverte $X_K\fleche X$,\cr
&X&\ktilde& $=X_k\times_k \ktilde$,&&\vtilde&&le prolongement de $v$ à
$\Ktildeetoile$,\cr
&\Xtilde&& $=X\times_\ogoth\ogothtilde$,&&I&&le groupe profini
$\Gal(\Kbar|\Ktilde)$,\cr
&S&&\multispan5\quad le $\Gal(\ktilde|k)$-ensemble des composantes
irréductibles de $X_{\ktilde}$.\hfil \cr
}
}

On suppose que $X_\Kbar$ est intègre.  Ce n'est qu'au \numero~5 que
nous supposons que $X_\Kbar$ est {\it birationnel\/} à $\P_2$ et que
l'application $\Pic X_\Ktilde\fleche\Pic X_\Kbar$ est un {\it
isomorphisme\/} --- s'y rapporter pour les hypothèses précises, qui
sont un peu plus générales.

\section{Localisation et filtrations}

Pour tout schéma $V$, désignons par $K_0(V)$ le groupe de Grothendieck
des $\Orond_V$-modules cohérents, par $(K_0(V)^{(t)})_{t\ge0}$ sa
filtration décroissante par la codimension du support, par $K^0(V)$
l'anneau de Grothendieck de $\Orond_V$-modules cohérents localement
libres muni de la $\gamma$-filtration $(K^0(V)^{(t)})_{t\ge0}$, par
$\theta : K^0(V)\fleche K_0(V)$ l'homomorphisme canonique de
$K^0(V)$-modules \citer\manin().  On a $\theta(K^0(V)^{(t)})\subset
K_0(V)^{(t)}$ pour $t\ge0$
\citer\sgavi(p.~523).

\th LEMME 1
\enonce
Les deux applications $\theta^{(1)}: K^0(X)^{(1)}\fleche K_0(X)^{(1)}$ 
et\/ $\theta^{(2)}: K^0(X)^{(2)}\fleche K_0(X)^{(2)}$ induites par
$\theta$ sont des isomorphismes.
\endth
\demonstration 
Comme $X$ est régulier et projectif sur $\ogoth$, l'homomorphisme
$\theta$ est bijectif \citer\manin().

Le groupe $K_0(X)^{(1)}$ est engendré par les faisceaux structuraux
$\Orond_V$ des fermés $V\subset X$ de codimension $\geq 1$ ; un
tel faisceau admet une résolution
$$
\zero \fleche \Erond_n \fleche \Erond_{n-1}
\fleche \cdots \fleche \Erond_1 \fleche
\Erond_0 \fleche \zero \numeroter
\newcount\resolution \global\resolution=\formuleno
$$
par des $\Orond_X$-modules $\Erond_t$ localement libres :
l'homologie de ce complexe s'annule en degré
$\neq0$ et s'identifie à $\Orond_V$ en degré $0$ ; l'image de l'élément
$x=\somme_{t=0}^n (-1)^t [\Erond_t]$
de $K^0(X)$ par $\theta$ vaut $[\Orond_V]\in K_0(X)$.  Comme le
complexe~\formule(\resolution) est exact au point générique de $X$,
on a $\rg(x)=0$.  Mais $K^0(X)^{(1)}$ est par définition le noyau de
l'homomorphisme $\rg:K^0(X)\rightarrow\Z$ ; on a donc
$x\in K^0(X)^{(1)}$ et l'on a établi la surjectivité de
$\theta^{(1)}$.

Pour ce qui concerne le cran 2 de la filtration, considérons le
composé
$$
\Pic X \droite{} {K^0(X)^{(1)}\over K^0(X)^{(2)}} \droite{}
{K_0(X)^{(1)}\over K_0(X)^{(2)}}
\numeroter
\newcount\compose \global\compose=\formuleno
$$
dans lequel la première flèche est l'isomorphisme
\citer\manin() qui envoie en
$[\Orond_X]-[\Lrond^{-1}]\in K^0(X)^{(1)}/K^0(X)^{(2)}$, quel que soit
le $\Orond_X$-module inversible $\Lrond$, la classe de $\Lrond$ dans
$\Pic X$ ; la seconde flèche est induite par ce qui précède.  On sait
que le composé \formule(\compose) est un isomorphisme
\citer\manin(p.~45) ; il en est donc de même de
$\theta^{(2)}:K^0(X)^{(2)}\fleche K_0(X)^{(2)}$, ce qu'il fallait
démontrer.\cqfd

L'immersion fermée $i:X_k\fleche X$ (resp. l'immersion ouverte~$j:
X_K\fleche X$) induit l'homomorphisme d'image directe $i_* :
K_0(X_k)\fleche K_0(X)$ (resp.~de restriction $j^* : K_0(X)\fleche
K_0(X_K)$) qui augmente la filtration par un cran (resp.~qui respecte
la filtration).

On note $v: \Ketoile \fleche \Z$ la valuation normalisée de $K$ et
$\pi$ une uniformisante de~$K$ ; on a $v(\pi)=1$.

\th LEMME 2
\enonce 
Les immersions $i$ et $j$ induisent une suite exacte 
$$
\advance\abovedisplayskip by 2pt 
K_0(X_k)^{(1)} \droite{i_*} K_0(X)^{(2)} \droite{j^*} 
K_0(X_K)^{(2)} \droite{} \zero.
\numeroter
$$
\endth
\newcount\filtree \filtree=\formuleno 
\demonstration 
Grâce à l'exactitude de la suite de localisation
\citer\borelserre(prop.~7) 
$$
\advance\abovedisplayskip by 2pt 
K_0(X_k)\fleche K_0(X) \fleche 
K_0(X_K)\fleche \zero, \numeroter
\newcount\localisation \global\localisation=\formuleno 
$$
il est clair que~\formule(\filtree) est un complexe et que $j^*$ est
surjectif ; il reste à faire voir que $i_*(K_0(X_k)^{(1)})$ contient
$\Ker(j^*)$. Soit $x\in\Ker(j^*)$ ; par l'exactitude
de~\formule(\localisation), il existe un $y\in K_0(X_k)$ tel que
$i_*(y)=x$.  On peut écrire $y=y_0+y_1$ où $y_0$ est l'image d'un
cycle $\xi$ de codimension~0 dans $X_k$ et où $y_1\in K_0(X_k)^{(1)}$
\citer\sgavi(p.~519). Comme $i_*(y_0)=x-i_*(y_1)$ appartient à
$K_0(X)^{(2)}$, l'image de $i_*(y_0)$ dans $\Pic(X)$ par
l'isomorphisme réciproque de~\formule(\compose) est nulle --- le
diviseur $\xi$ de $X$ est principal.

Écrivons alors $\xi=\div_X(f)$, où $f\in K(X)^\times$ ; on a
$\div_{X_K}(f)=0$ et, comme $X_\Kbar$ est projectif et intègre,
$f\in \Ketoile$.  Ainsi, $\xi=v(f)\div_X(\pi)$.  En prenant les
classes dans $K_0(X_k)$, on a $y_0=v(f)[\Orond_{X_k}]+z_1$ pour
$z_1\in K_0(X_k)^{(1)}$ \citer\sgavi(p.~520, 1.1.3). Ainsi
$x=v(f)[\Orond_{X_k}]+i_*(y_1+z_1)$.  Mais on a $[\Orond_{X_k}]=0$
dans $K_0(X)$ à cause de l'exactitude de la suite
$$
\zero\fleche\Orond_X \droite{\pi}\Orond_X \fleche \Orond_{X_k}\fleche\zero.
$$
Ceci montre que $\Ker(j^*)\subset\Im(i_*)$ et
établit par là l'exactitude de~\formule(\filtree).\cqfd 
\goodbreak

\section{L'homomorphisme de spécialisation}

Soit $i_Y:Y\fleche X$ un diviseur effectif à support dans $X_k$.  On
désigne par $i_Y^{\star}$ le composé
$$
\advance\abovedisplayskip by 2pt 
K_0(X) \hfl{\theta^{-1}}{\hbox{\~{}}}{8mm} K^0(X) \droite{i_Y^*} K^0(Y)
\droite{\theta} K_0(Y) \numeroter
\newcount\specialisation  \global\specialisation=\formuleno
$$
(« homomorphisme de spécialisation »).  La flèche composée
$i_Y^{\star}$~\formule(\specialisation) de spécialisation est à
distinguer de la flèche médiane $i_Y^*: K^0(X)\fleche K^0(Y)$.
Lorsque $Y=X_k$, on abrège $i_{X_k}^{\star}$ simplement en
$i^{\star}$.

\th LEMME 3
\enonce  
Soit $\Frond$ un $\Orond_X$-module cohérent.  Dans le groupe $K_0(Y)$,
on a
$$
i_Y^{\star}([\Frond]) = \somme_{t\geq 0}(-1)^t 
[\tor_t^{\Orond_X}(\Frond,\Orond_Y)]. \numeroter
$$\vskip-\lastskip
\endth
\newcount\torsalt \torsalt=\formuleno 
\demonstration  Soit \formule(\resolution) une résolution de $\Frond$
par des $\Orond_X$-modules $\Erond_t$ localement libres ; on a 
$i_Y^{\star}([\Frond])=\somme_{t=0}^n (-1)^t [\Erond_t\otimes \Orond_Y]$, 
ce qui implique~\formule(\torsalt) car
les $\tor_t^{\Orond_X}(\Frond,\Orond_Y)$ sont les groupes
d'homologie du complexe $(\Erond_t\otimes \Orond_Y)_t$.\cqfd

\goodbreak

\th LEMME 4
\enonce  
Notons $\Irond$ l'idéal inversible définissant $Y$ dans $X$ et
$\Nrond=\Irond/\Irond^2$ son faisceau conormal.  Alors
$$
i_Y^{\star} i_{Y*}^{\phantom{\star}}(x) = ([\Orond_Y] -[\Nrond]).x
\qquad (x\in K_0(Y)).
\numeroter
$$\vskip-\lastskip
\endth
\newcount\conormal \global\conormal=\formuleno 
\demonstration D'après \formule(\torsalt) et par
linéarité, il suffit de faire voir que pour tout $\Orond_Y$-module
cohérent $\Frond$, on a  
$$
\tor_1^{\Orond_X}(\Frond,\Orond_Y) =\Frond \otimes_{\Orond_X}
\Nrond\quad\hbox{et}\quad 
\tor_t^{\Orond_X}(\Frond,\Orond_Y) =\zero\quad \hbox{pour}\ t\geq2.
\numeroter \newcount\torsup \global\torsup=\formuleno
$$
À la suite exacte tautologique
$\zero\fleche\Irond\fleche\Orond_X\fleche\Orond_Y\fleche\zero$ est
associée la suite exacte bornée à droite
$$
\cdots\fleche \tor_t^{\Orond_X}(\Frond,\Orond_X)\fleche
\tor_t^{\Orond_X} (\Frond,\Orond_Y)\fleche
\tor_{t-1}^{\Orond_X}(\Frond,\Irond)\fleche\cdots.
\numeroter\newcount\longuetor \global\longuetor=\formuleno
$$
Le $\Orond_X$-module $\Irond$ étant inversible,
$\tor_t^{\Orond_X}(\Frond,\Orond_X)=\zero$ et
$\tor_t^{\Orond_X}(\Frond,\Irond)=\zero$ pour $t\geq 1$, ce qui donne
la seconde égalité~\formule(\torsup). D'autre part, tensorisant la
suite tautologique avec $\Irond$, on obtient l'isomorphisme
$\Nrond\fleche\Irond\otimes_{\Orond_X}\Orond_Y$ qui, joint à
l'isomorphisme
$\tor_1^{\Orond_X}(\Frond,\Orond_Y)\fleche\Frond\otimes_{\Orond_X}\Irond$
dans \formule(\longuetor) et au fait que le $\Orond_X$-module
${\Frond}$ provient de $Y$, implique la première
égalité~\formule(\torsup). \cqfd

\goodbreak

\th LEMME 5
\enonce  On a $i^{\star}\circ i_*=0$. Il existe un homomorphisme et un seul 
$\sigma : K_0(X_K) \fleche K_0(X_k)$ tel que $i^{\star}=\sigma\circ j^*$.
\endth
\demonstration  Appliquons \formule(\conormal) avec $Y=X_k$ ; comme 
le $\Orond_X$-module inversible $\Irond=\pi\Orond_X$ est alors libre,
il en de même du $\Orond_Y$-module inversible
$\Nrond=\Irond\otimes_{\Orond_X}\Orond_Y$ et par conséquent
$[\Orond_Y]-[\Nrond]=0$.  L'existence et l'unicité de $\sigma$ en
résultent, compte tenu de l'exactitude de
\formule(\localisation). \cqfd

\goodbreak

\th LEMME 6
\enonce  
La spécialisation $i_Y^{\star}=\theta\circ i_Y^*\circ\theta^{-1}$
\formule(\specialisation) respecte la filtration.
\endth
\demonstration   
Clairement ${i_Y^*}$ respecte le cran~1 de la filtration puisqu'il
préserve les rangs ; il respecte le cran~2 à cause de la
commutativité du diagramme
$$
\def\droite#1{\!\!\hfl{#1}{}{12mm}\!\!}
\diagram{
\zero\fleche&K^0(X)^{(2)}
&\fleche&K^0(X)^{(1)}&\droite{\det_X}&\Pic(X)&\fleche\zero\cr  
&&&\versbas{i_Y^*}&&\versbas{i_Y^*}\cr
\zero\fleche&K^0(Y)^{(2)}
&\fleche&K^0(Y)^{(1)}&\droite{\det_Y}&\Pic(Y)&\fleche\zero\cr}
$$
dont les lignes sont exactes puisque $Y$ est localement d'intersection
complète dans $X$ régulier \citer\manin(p.~45).  Vu le lemme~1, ceci
implique que $i_Y^{\star}$ respecte les crans~1 et~2 de la filtration.

Reste à faire voir que $i_Y^{\star}(K_0(X)^{(3)})=0$.  Soient $P\in X$
un point de codimension~3 et $x\in K_0(X)^{(3)}$ sa classe.  Si
$P\notin Y$, on a $i_Y^{\star}(x)=0$.  Si $P\in Y$, soit $y\in K_0(Y)$
sa classe.  D'après~\formule(\conormal), on a
$i_Y^{\star}(x)=([\Orond_Y]-[\Nrond]).y$ où $\Nrond$ est le faisceau
conormal de $Y$ dans $X$.  Or $[\Orond_Y] -[\Nrond]$ appartient à
$K^0(Y)^{(1)}$ car son rang est nul et $y$ appartient à $K_0(Y)^{(2)}$
donc $i_Y^{\star}(x)$ appartient à $K_0(Y)^{(3)}=\zero$
\citer\sgavi(p.~523), ce qu'il fallait démontrer.\cqfd

\goodbreak

Soit $Y$ une composante irréductible de $X_k$.  Notons $\check{Y}$ le
schéma réduit associé à $Y$ et $m_Y$ la multiplicité de $\check{Y}$
dans $X_k$.  Nous allons considérer $Y$ comme un {\it sous-schéma
fermé\/} de $X_k$ --- définie par le faisceau d'idéaux $\Irond^{m_Y}$, où
$\Irond$ est le faisceau d'idéaux définissant $\check{Y}$ dans $X_k$.
Notons $r_Y$, $g_Y$, $f_Y$ les morphismes canoniques suivants : 
$$
r_Y : \check Y\fleche Y,\ \ \ 
g_Y : \check Y\fleche X_k,\ \ \ 
f_Y : Y\fleche X_k.
\numeroter
\newcount\rygyfy  \global\rygyfy=\formuleno
$$

\th LEMME 7
\enonce
On a $\;i_Y^{\star}=(m_Yr_{Y*})\circ
i_{\check{Y}}^{\star}\;:\;K_0(X)^{(2)}\fleche K_0(Y)^{(2)}$. 
\endth
\demonstration   
Écrivons la formule de projection pour
$r_Y:\check Y\fleche Y$~\formule(\rygyfy) 
$$
r_{Y*}(r_Y^*(y). x)=y. r_{Y*}(x)\qquad (x\in K_0(\check Y),\ y \in
K^0(Y)) \numeroter\newcount\projection \global\projection=\formuleno
$$
\citer\sgavi(p.~287).  Prenons $x=[\Orond_{\check Y}]$ et
$y=i_Y^*\circ \theta^{-1}(z)$ avec 
$z\in K_0(X)$.  On a alors $r_Y^*(y).x=i_{\check{Y}}^{\star}(z)$ et
$\theta(y)=y. [\Orond_Y]=i_{Y}^{\star}(z)$.
Multiplions~\formule(\projection) par~$m_Y:$
$$
m_Yr_{Y*}(i_{\check{Y}}^{\star}(z))=i_{Y}^{\star}(z)-
y. ([\Orond_Y]-m_Yr_{Y*}[\Orond_{\check Y}]). 
$$
Si $z\in K_0(X)^{(2)}$, on a $y\in K^0(Y)^{(2)}$ (lemme~1) et la
définition de la multiplicité $m_Y$ implique que
$[\Orond_Y]-m_Yr_{Y*}[\Orond_{\check Y}] \in K_0(Y)^{(1)}$
\citer\sgavi(p.~519), donc
$y.([\Orond_Y]-m_Yr_{Y*}[\Orond_{\check Y}])\in K_0(Y)^{(3)}=\zero$
\citer\sgavi(p.~523), d'où le lemme. \cqfd

\goodbreak

\th LEMME 8
\enonce  
L'application $i^\star:K_0(X)^{(2)}\fleche K_0(X_k)^{(2)}$ (lemme\/~$6$)
coïncide avec le composé
$$
\advance\abovedisplayskip by 4pt
\def\droite#1{\;\hfl{#1}{}{15mm}\;}
K_0(X)^{(2)}\droite{(i_Y^\star)_Y^{\phantom{\star}}} 
\bigoplus_Y K_0(Y)^{(2)}
\droite{\sum_Y f_{Y*}} K_0(X_k)^{(2)}
\advance\belowdisplayskip by -2pt
$$
et aussi avec le composé\/
$\left(\sum_Ym_Yg_{Y*}\right)\circ(i_{\check Y}^\star)_Y^{\phantom{\star}}$,
où $Y$ parcourt les composantes irréductibles de $X_k$.
\endth
\demonstration   
Elle est identique à celle du lemme~7, à cela près qu'ici on se sert
du fait que
$[\Orond_{X_k}]-\sum_Yf_{Y_*}([\Orond_Y])\in K_0(X_k)^{(1)}$
(cf.~\citer\sgavi(p.~519)). \cqfd

Pour tout diviseur effectif\/ $Y$ de $X$ à support dans $X_k$, on
note $\chi_{Y} : K_0(Y)^{(2)} \fleche \Z$ la caractéristique
d'Euler-Poincaré.  Lorsque $Y=X_k$, on abrège $\chi_{X_k}$
simplement en $\chi$.

\th LEMME 9
\enonce  
Pour tout $x\in K_0(X)^{(2)}$ (resp.~tout $y\in K_0(X_k)^{(1)}$),  on a
$$
\chi(i^\star(x)) =\somme_Y m_Y \chi_{\check Y} (i_{\check
Y}^\star(x))\quad\left(\hbox{resp.\ }
\somme_Y m_Y
\chi_{\check Y}(i_{\check Y}^\star \circ i_*(y))=0\right)
$$
dans $\Z$, où $Y$ parcourt les composantes
irréductibles de $X_k$.   
\endth
\demonstration
Vu que
$\chi\left(\somme_Ym_Yg_{Y*}\right)=\somme_Ym_Y\chi_{\check Y}$, la
première assertion résulte du lemme~8.  La deuxième suit en prenant
$x=i_*(y)$, vu le lemme~5.
\cqfd

Puisque l'extension $\Ktilde|K$ est non ramifiée, le schéma $\Xtilde$
est régulier et {\it tout ce qui précède s'applique à\/ $\Xtilde$
sur\/ $\ogothtilde$} --- ce fait sera utilisé par la suite sans y
faire allusion.

Désignons par $(e_Z)_{Z\in S}$ la base canonique du $\Z$-module libre
$\Z^S$, où $S$ est l'ensemble des composantes irréductibles rendues
réduites de $\Xktilde$.  Soit $P$ l'élément $\somme_{Z\in S}m_Z{e_Z}$
(où $m_Z$ est la multiplicité de $Z$ dans $\Xktilde$) de $\Z^S$ et
soit $\xi:\Hom_k(\Z^S,\Z)\fleche \Z$ l'homomorphisme
$h\longmapsto h(P)$.

Notant $x\mapsto\xtilde$ l'application
$K_0(X)^{(2)}\fleche K_0(\Xtilde)^{(2)}$, on dispose des
homomorphismes
$$
\eqalign{
\psitilde&:K_0(\Xtilde)^{(2)}\fleche\Hom(\Z^S,\Z),&
\psitilde(x)(e_Z)=\chi_Z(i_Z^\star(x))\cr
\psi&:K_0(X)^{(2)}\fleche\Hom_k(\Z^S,\Z),&
\psi(x)(e_Z)=\chi_Z(i_Z^\star(\xtilde))\cr
}\numeroter\newcount\psidef \global\psidef=\formuleno
$$

Le quotient $B(X)$ de $\Hom_k(\Z^S,\Z)$ par
$\psi(i_*(K_0(X_k)^{(1)}))$ joue un rôle important dans la suite.

\th LEMME 10
\enonce
Pour tout $x\in K_0(X_k)^{(1)}$, on a $\;\xi(\psi(i_*(x)))=0.$
\endth
\demonstration
Il suffit d'appliquer le lemme~9 à $\Xtilde$.\cqfd 

D'après le lemme~10, il y a un homomorphisme et un seul
$\xibar:B(X)\fleche\Z$ qui induit $\xi$ sur $\Hom_k(\Z^S,\Z)$;
nous allons désigner le noyau de $\xibar$ par $B(X)_0$.

On note $T$ le quotient de $\Z^S$ par le sous-groupe engendré par
$P$ ; il est muni d'une action de $\Gal(\ktilde|k)$ ; $\Hom_k(T,\Z)$
s'identifie à un sous-groupe de $\Hom_k(\Z^S,\Z)$.

\th LEMME 11
\enonce
On a $\psi(i_*(K_0(X_k)^{(1)}))\subset\Hom_k(T,\Z)$ ; le
quotient est égal à
$B(X)_0=\Ker(\xibar:B(X)\fleche\Z)$.
\endth
\demonstration
Vu le lemme~10, cela résulte de ce que $\Hom_k(T,\Z)$
est le noyau de $\xi :\Hom_k(\Z^S,\Z)\fleche\Z$. \cqfd

\section{La surjectivité de la spécialisation}

Soit $Y$ une composante irréductible de $X_k$ et soit
$\Kprim\subset\Ktilde$ une extension finie de $K$ telle que toutes les
composantes irréductibles de $Y_\ktilde$ soient définissables sur
$\kprim$, le corps résiduel de $\Kprim$.  On note $\ogothprim$
l'anneau des entiers de $\Kprim$ et l'on pose
$\Xprim=X\times_\ogoth\ogothprim$.

Soit $Z$ une composante irréductible (rendue réduite) de $Y_\kprim$ et
soit $m$ sa multiplicité.  Les deux lemmes suivants sont adaptés de
\citer\neron(p.~240--242).

\th LEMME 12
\enonce
Il existe un ouvert dense\/ $U\subset Z$ tel que pour toute extension
finie\/ $l$ de\/ $\kprim$ et pour tout point\/
$x\in U(l)$, l'anneau local\/ $\Orond_{Z_l,x}$ admet un système de
paramètres $(f_1,f_2)$ vérifiant\/ $\dim_l \Orond_{Z_l,x}/(f)=m$.
\endth
\demonstration Partons d'un ouvert dense $U\subset Z$ tel que le
schéma réduit $\check U$ soit lisse sur $\kprim$, de sorte que le
$\Orond_{{\check U}}$-module $\Omega_{{\check U}|\kprim}$ est
localement libre.  Remplaçant $U$ par un ouvert dense, on peut
supposer que le $\Orond_{{\check U}}$-module
$\Omega_{{\check U}|\kprim}$ est libre.  Comme
$\Omega_{{\check U}|\kprim}$ est un quotient de $\Omega_{{U}|\kprim}$,
on peut trouver $f_1,f_2\in \Gamma(U,\Orond_U)$ tels que les images
des différentiels $df_1,df_2$ dans $\Omega_{{\check U}|\kprim}$
forment une base de ce module.  La restriction du morphisme
$$
f=(f_1,f_2):U\fleche V,\quad V={\bf A}^2_{\kprim}
$$
au schéma réduit ${\check U}$ étant étale, on peut supposer que $f$
est fini et plat, quitte à remplacer $U$ et $V$ par des ouverts
denses.  Montrons que $U$ convient.  Soient $l$ une extension finie et
$x\in U(l)$ un point $l$-rationnel; on peut supposer que $f_l(x)$ est
l'origine de~$V_l=V\times_{\kprim}l$.  Soit ${}^hV_l$ le hensélisé de
$V_l$ à l'origine et soit ${}^hU_l$ la composante locale de
$U_l\times_{V_l} {}^hV_l$ au-dessus de~$x$.  Le morphisme induit
${}^h\!{\check f_l}:{}^h{\check U_l}\fleche {}^hV_l$ étant un
isomorphisme, le degré de ${}^h\!f_l$ vaut la longueur~$m$ de l'anneau
local du point générique de ${}^hU_l$, ce qu'il fallait
démontrer. \cqfd
\goodbreak

Nous utilisons la variante suivante du lemme classique de Hensel : 

\th LEMME 13
\enonce
Il y a un ouvert dense\/ $U\subset Z$ tel que pour toute extension
finie\/ $l\subset\ktilde$ de\/ $\kprim$, tout point\/ $l$-rationnel\/
$x\in U_l(l)$ se relève : il existe un\/ $\Ogoth$-schéma\/ $C$ plat et
fini de rang\/~$m$ et une\/ $\Ogoth$-immersion fermée\/
$a:C\fleche X_\Ogoth$ tels que\/ $a_l(C_l)=x$, où\/ $\Ogoth$ est
l'anneau des entiers de l'extension\/ $L\subset\Ktilde$ de corps
résiduel\/ $l$.
\endth
\demonstration  Soit $W$ un ouvert de $\Xprim$ tel que $W_\kprim$ soit
un ouvert dense de $Z$. Nous dirons d'un énoncé qu'il est vrai « au
rétrécissement de $W$ près » pour signifier qu'il est vrai en
remplaçant $W$ par un ouvert $O\subset W$ tel que $O_\kprim$ soit
dense dans $Z$.  Au rétrécissement de $W$ près, pour toute extension
finie\/ $l\subset\ktilde$ de\/ $\kprim$ et pour tout point\/
$x\in U(l)$, l'anneau local $\Orond_{Z_l,x}$ admet un système de
paramètres $\overline f$ tel que
$\dim_l \Orond_{Z_l,x}/({\overline f})=m$ (lemme~12).  Au
rétrécissement de $W$ près, ${\overline f}$ se relève en une suite
$f$ d'éléments de $\Gamma(W_\Ogoth,\Orond_{W_\Ogoth})$ ; c'est une
suite régulière de $\Orond_{W_\Ogoth,x}$ \citer\ega(0$_{IV}$,
15.1.16).  Au rétrécissement de $W$ près, une composante locale $C$
de $V(f)$ contenant~$x$ est finie et plate sur $\Ogoth$ ; ce $C$
convient \citer\ega(0$_{IV}$, 15.1.16).  \cqfd

\th LEMME 14
\enonce
Il y a un ouvert dense $U\subset Z$ tel que tout point fermé $x\in U$
se relève : il existe un $\ogothprim$-schéma $C$ plat et fini de rang
$m\deg(x)$ et une $\ogothprim$-immersion fermée $C\fleche\Xprim$ tel
que $C_\kprim=x$.
\endth
\demonstration  Montrons qu'un ouvert $U\subset Z$ fourni par le
lemme~13 convient.  Soit $x\in U$ un point fermé.  Identifions le
corps résiduel $\kprim(x)$ à $l\subset\ktilde$ et choisissons un point
$l$-rationnel $y$ de $U_l$ au-dessus de $x$.  Il existe un
$\Ogoth$-schéma (où $\Ogoth$ est l'anneau des entiers de l'extension
$L\subset\Ktilde$ de corps résiduel $l$) plat et fini de rang $m$ et
une $\Ogoth$-immersion fermée $a:C\fleche X_\Ogoth$ tels que
$a_l(C_l)=y$ (lemme~13).  Or $C$ est un $\ogothprim$-schéma fini et
plat de rang $m\deg(x)$ et le composé $C\fleche X_\Ogoth\fleche\Xprim$
une $\ogothprim$-immersion $C\fleche\Xprim$ tels que
$C_\kprim=x$. \cqfd

\th PROPOSITION 1
\enonce
L'application $\psi:K_0(X)^{(2)}\fleche\Hom_k(\Z^S,\Z)$
\formule(\psidef) est surjective.
\endth
\demonstration  Soit $(b_Y)_{Y}$ la base canonique de
$\Hom_k(\Z^S,\Z)$, indexée par les composantes irréductibles de
$X_k$ : pour $Z\in S$, on a $b_Y(e_Z)=1$ si $Z$ est une composante
irréductible de $Y_\ktilde$ et $b_Y(e_Z)=0$ sinon.  Soit $Y$ une
composante irréductible de $X_k$ et soit $Z$ une composante
irréductible de $Y_\kprim$ ($\kprim|k$ est une extension finie ---
corps résiduel de $\Kprim\subset\Ktilde$ --- où toutes les composantes
irréductible de $Y_\ktilde$ sont définissables) ; soit
$H\subset\Gal(\kprim|k)$ le stabilisateur de $Z$.  Tout ouvert dense
de $Z$ possède un 0-cycle $H$-stable de degré~1 (Lang-Weil);
prenons-en un qui se relève en un 1-cycle $H$-stable $C$ de
$\Xprim=X\times_\ogoth\ogothprim$, où $\ogothprim$ est l'anneau des
entiers de $\Kprim$ (lemme 14).  Pour un conjugué $^\tau Z$
($\tau\in\Gal(\kprim|k)/H)$ de $Z$, posons $C_{\tau}={}^\tau C$.
Alors
$$
C_Y=\somme_\tau  C_{\tau},\quad\tau\in\Gal(\kprim|k)/H
$$
est un 1-cycle de $\Xprim$ stable par $\Gal(\kprim|k)$ et qui provient
donc d'un 1-cycle de $X$; soit $y\in K_0(X)^{(2)}$ sa classe. Comme
$\psi(y)=b_Y$, on a gagné.\cqfd

Le groupe $K_0(X_K)^{(2)}$ s'identifie à $\azero{X_K}$ et le noyau de
$\chi_{X_K}$ au noyau $\azeroo{X_K}$ de l'application degré
$\deg:\azero{X_K}\fleche\Z$. 

\goodbreak

\th PROPOSITION 2
\enonce
Avec ces notations, il existe un unique homomorphisme
$$
\gamma:\azero{X_K}\droite{} B(X) 
= {\Hom_k(\Z^S,\Z)\over\psi(i_*(K_0(X_k)^{(1)}))}
\numeroter
$$
tel que, pour tout\/ $y\in K_0(X)^{(2)}$, la classe de\/
$\psi(y)\in\Hom_k(\Z^S,\Z)$~\formule(\psidef) dans\/ $B(X)$ soit égale
à\/ $\gamma(j^*(y))$; cet homomorphisme est surjectif.
\endth
\newcount\gammadef \global\gammadef=\formuleno
\demonstration L'existence et l'unicité de $\gamma$ résultent aussitôt
de la définition de $\psi$~\formule(\psidef) et de l'exactitude de la
suite~\formule(\filtree); la surjectivité de $\gamma$ résulte de la
prop.~1.\cqfd

\th LEMME 15
\enonce
Pour tout\/ $x\in\azero{X_K}$, on a\/
$\deg(x)=\xibar(\gamma(x))$.
\endth
\demonstration  Soit $y\in K(X)^{(2)}$ un relèvement de~$x$.  Nous avons
$$
\eqalign{
\xibar(\gamma(x))&=\psi(y)(P)\cr
&=\somme_{Z\in S} m_Z\chi_Z(i_Z^\star(\tilde y))\cr
&=\chi_{X_k} (i^\star(y))\quad\hbox{(lemme 9)}\cr
&=\chi_{X_K} (x).\quad\citer\sgavi(\hbox{p.~563})\qquad\qed\cr}
$$

Grâce au lemme~15 (cf.~lemme~11), $\gamma$ définit un homomorphisme
$$
\gamma_0:\azeroo{X_K}\droite{} B(X)_0 =
{\Hom_k(T,\Z)\over\psi(i_*(K_0(X_k)^{(1)}))}.
\numeroter
\newcount\gammaodef \global\gammaodef=\formuleno
$$

\th PROPOSITION 3
\enonce
L'application $\gamma_0$ \formule(\gammaodef) est surjective. 
\endth
\demonstration C'est une conséquence immédiate de la prop.~2.\cqfd

\section{Comparaison avec l'application caractéristique}

Nous allons noter $\Zrond^t(V)$ le groupe des cycles de codimension
$t$ sur un schéma $V$.  La valuation canonique de $\Ktilde$ est notée
$\vtilde$.

\th LEMME 16
\enonce 
Le noyau du composé\/ $\Z^S\fleche\Zrond^1(\Xtilde)\fleche\Pic\Xtilde$
est engendré par\/ $P=\div_\Xtilde(\pi)$ ; son conoyau s'identifie à\/
$\Pic X_\Ktilde$ : on a la suite exacte
$$
\zero\fleche T\fleche \Pic \Xtilde \droite{j^*} \Pic X_\Ktilde\fleche \zero.
\numeroter
$$
\endth
\newcount\suitepic \global\suitepic=\formuleno
\demonstration C'est une petite chasse au diagramme.  On se sert du
fait que les seules fonctions inversibles sur $\Xtilde$
(resp.~$\XKtilde$) sont $\ogothtildemult$
(resp.~$\Ktildeetoile$).\cqfd

À la suite exacte \formule(\suitepic) de $\Gal(\ktilde|k)$-modules
est associée la suite exacte 
$$
\advance\abovedisplayskip by 4pt
\Hom_k(\Pic \Xtilde,\Z)\fleche\Hom_k(T,\Z)
\droite{\delta}\Ext_k^1 (\Pic X_\Ktilde,\Z).
\numeroter\newcount\defdelta \global\defdelta=\formuleno
$$

\th LEMME 17
\enonce 
Pour tout $\,x\in K_0(X_k)^{(1)}$, l'on a $\;\delta(\psi(i_*(x)))=0$ : il
existe un unique homomorphisme
$\;\deltabar:B(X)_0\fleche\Ext_k^1 (\Pic X_\Ktilde,\Z)$
qui induit $\;\delta$~\formule(\defdelta) sur\/ $\Hom_k(T,\Z)$.
\endth
\demonstration 
D'après l'exactitude de la suite~\formule(\defdelta), il suffit
de montrer la commutativité du carré
$$
\diagram{
\Zrond^1(X_k) &\hskip-8mm \hfl{}{}{26mm} &\hskip-2mm K_0(X_k)^{(1)}\cr
\vfl{}{}{10mm}&&\hskip-2mm \vfl{}{\psi\circ i_*}{10mm}\cr
\Hom_k(\Pic \Xtilde,\Z)&\hskip-2mm\hfl{}{}{20mm}&\hskip-2mm\Hom_k(T,\Z)\cr
}
\numeroter\newcount\lecarre \global\lecarre=\formuleno
$$
dans lequel la verticale de gauche est l'homomorphisme qui à une
courbe irréductible $C\in X_k^{(1)}$ tracée sur $X_k$ 
fait correspondre le composé
$$
\advance\abovedisplayskip by 5pt 
\Pic \Xtilde\hfl{}{}{10mm}\Pic \Ctilde\hfl{\deg_{\Ctilde}}{}{12mm}\Z, 
\quad\hbox{où}\ \Ctilde=C\times_k \ktilde
$$
et où $\deg_{\Ctilde}$ désigne le « degré total ».  Montrons plus
généralement que pour une courbe irréductible $D\in X_\ktilde^{(1)}$
(de classe $[\Orond_D]\in K_0(\Xtilde)$) et une composante
irréductible $Z\in X_\ktilde^{(0)}$, on a
$$
\chi_Z i_Z^\star ([\Orond_D])=\deg_D(\Orond_\Xtilde(Z)|_D).
\numeroter\newcount\degdeg \global\degdeg=\formuleno
$$

Tensorisant la suite exacte $\zero \fleche \Orond_{\Xtilde}(-Z) \fleche \Orond_{\Xtilde} \fleche\Orond_Z \fleche \zero$ de $\Orond_\Xtilde$-modules cohérents avec $\Orond_D$, on
obtient la suite exacte 
$$
\zero\fleche \toroxtilde_1(\Orond_D,\Orond_Z)  \fleche  
\Orond_\Xtilde(-Z)|_D\fleche\Orond_D  \fleche
\Orond_D\otimes_{\Orond_\Xtilde} \Orond_Z \fleche  \zero  
\numeroter\newcount\longtor \global\longtor=\formuleno
$$
dans laquelle les $\toroxtilde_t(\Orond_C,\Orond_Z)$
s'annulent pour $t\geq 2$ parce que les $\Orond_\Xtilde$-modules
$\Orond_\Xtilde$ et $\Orond_\Xtilde(-Z)$ sont localement libres.
Il vient
$$
\eqalign{
\chi_Z i_Z^\star ([\Orond_D])&=\chi_Z(\Orond_D\otimes_{\Orond_\Xtilde}\Orond_Z)
-\chi_Z(\toroxtilde_1(\Orond_D,\Orond_Z))\qquad\hbox{(lemme 3)}\cr
&=\chi_D(\Orond_D\otimes_{\Orond_\Xtilde}\Orond_Z)
-\chi_D(\toroxtilde_1(\Orond_D,\Orond_Z))\cr 
&=\chi_D(\Orond_D)-\chi_D(\Orond_\Xtilde(-Z)|_D)
\qquad\hbox{d'après}\ \formule(\longtor) \cr
&=\deg_D(\Orond_\Xtilde(Z)|_D)\qquad\hbox{(cf.~\citer\neron(p.~238, th.~1))}\cr
}
$$
ce qui établit~\formule(\degdeg) et par là la commutativité
de~\formule(\lecarre).\cqfd

Soient $\Spec L\fleche \XKtilde$ un point fermé de $\XKtilde$ (d'image
$Q\in\XKtilde$), 
$u:\Spec\Ogoth\fleche \Xtilde$ (d'image $C$, l'adhérence schématique
de $Q$ dans $\Xtilde$) son prolongement au spectre de l'anneau des
entiers $\Ogoth$ de $L$ ; le corps de fonctions sur~$C$ s'identifie
à~$L$.  Soient $Z$ une composante irréductible (rendue réduite) de
$\Xktilde$ et $f$ une équation locale du diviseur $Z$ de~$\Xtilde$ sur
un ouvert affine $\Spec A\subset\Xtilde$ contenant~$C$.

\th LEMME 18
\enonce 
On a $\;\chi_Zi_Z^\star([\Orond_C])=w(u^*(f))$ où\/
$w:L^\times\fleche\Z$ est la valuation normalisée de\/ $L$.
\endth
\demonstration   L'image
$u^*(f)$ de $f\in A$ dans l'anneau quotient $B=\Gamma(C,\Orond_C)$
correspondant \`a $C$ n'est pas diviseur de z\'ero puisque ce dernier
est int\`egre, vu qu'il poss\`ede le point g\'en\'erique $u(\Spec L)$,
et puisque $u^*(f)$ est non nul, $C$ n'\'etant pas contenu dans $Z$.
Il en résulte que $\toroxtilde_t(\Orond_C,\Orond_Z)=0$ pour $t\ge1$
(cf.~\citer\manin(p.~9, th.~2.3)), et le lemme~3 implique alors que
$$
\chi_Z i_Z^\star([\Orond_C])
=\chi_Z([\Orond_C \otimes_{\Orond_\Xtilde}\Orond_Z]).
$$
Le faisceau $\Orond_C \otimes_{\Orond_\Xtilde}\Orond_Z$
ayant son support dans $C\cap Z$ et ce dernier \'etant de dimension
$\le0$, on a $H^t(Z,\Orond_C\otimes_{\Orond_\Xtilde}\Orond_Z)=0$ pour
$t\ge1$, de sorte que 
$$
\eqalign{
\chi_Z(\Orond_C\otimes_{\Orond_\Xtilde}\Orond_Z)
&=\long H^0(Z,\Orond_C \otimes_{\Orond_\Xtilde}\Orond_Z)\cr
&=\long B/u^*(f)B.
}
$$
D'autre part, on a $w(u^*(f))=\long \Ogoth/u^*(f)\Ogoth$ et il reste à
faire voir que $B/u^*(f)B$ et $\Ogoth/u^*(f)\Ogoth$ ont la même
longueur.  Cela résulte aussitôt du diagramme commutatif de
$B$-modules
$$
\diagram{
\zero &\hfl{}{}{6mm} &B &\hfl{}{}{6mm} &\Ogoth &\hfl{}{}{6mm} &R
&\hfl{}{}{6mm} &\zero\cr
&&\vfl{}{u^*(f)}{4mm} &&\vfl{}{u^*(f)}{4mm} &&\vfl{}{\overline
f}{4mm}\cr
\zero &\hfl{}{}{6mm} &B &\hfl{}{}{6mm} &\Ogoth &\hfl{}{}{6mm} &R
&\hfl{}{}{6mm} &\zero\cr
}
$$
($R$ étant le $B$-module quotient $\Ogoth/B$) qui donne la suite
exacte 
$$
\zero \fleche {}_{\overline f}R \fleche
B/u^*(f)B\fleche\Ogoth/u^*(f)\Ogoth \fleche 
R/\overline f R \fleche \zero 
$$
$(u^*(f)$ n'\'etant point diviseur de z\'ero) : on a
$\long{}_{\overline f}R=\long R/\overline f R$ par 
l'exactitude de la suite
$$
\zero \fleche {}_{\overline f}R \fleche R \droite{\overline f}
R \fleche R/\overline f R \fleche \zero
$$
donc  $\long B/u^*(f)B =\long\Ogoth/u^*(f)\Ogoth$, démontrant le lemme.\cqfd

Notons $\Atilde$ l'anneau local de $\XKtilde$ le long de $Q$,
$\alpha:\Atildeetoile\fleche \Z^S$ l'homomorphisme
$g\mapsto\pr\circ\div_{\Xtilde}(g)$ où $\pr$ est la projection
$\Zrond^1(\Xtilde)\fleche \Z^S$, $N:\Letoile\fleche\Ktildeetoile$
l'homomorphisme « norme » et $\ev:\Atildeetoile\fleche\Ktildeetoile$
l'application d'évaluation en~$Q$ (on a $\;\ev(g)=N(u^*(g))$ pour tout
$g\in\Atildeetoile$).

\th LEMME 19
\enonce
On a $\;\vtilde\circ\ev=\psitilde([\Orond_C])\circ\alpha$ en tant
qu'applications $\;\Atildeetoile\fleche\Z$.
\endth
\demonstration Il suffit de considérer le cas
d'une équation locale $f\in\Atildeetoile$ d'une composante
irréductible $Z$ de $\Xktilde$ ($\alpha(f)=e_Z)$). Il vient alors
$$
\eqalign{
\vtilde\circ\ev(f)&=\vtilde\circ N\circ u^*(f)&(\ev=N\circ u^*)\cr
&=w(u^*(f))&(w=\vtilde\circ N)\cr
&=\chi_Zi_Z^\star([\Orond_C])&(\hbox{lemme 18})\cr
&=\psitilde([\Orond_C])(e_Z)&\hbox{d'après\ }\formule(\psidef)\cr
&=\psitilde([\Orond_C])\circ\alpha(f).&\hbox{\qed}\cr
}
$$

Soit $c$ un \ocycle\ sur $X_\Ktilde$, $\Atilde$ l'anneau semilocal de
$X_\Ktilde$ le long de $c$ et
$\ev:\Atilde{}^\times\fleche \Ktildeetoile$ l'application
d'évaluation en~$c$ $\bigl($en écrivant $c=\somme_Q n_Q Q$ où $Q$
parcourt les points fermés de $X_\Ktilde$ et où les $n_Q\in\Z$ sont
presque tous nuls, on a $\ev(g)=\prod\limits_Q\ev(Q)^{n_Q}\bigr)$.
Soit $C_Q$ l'adhérence schématique de $Q$ dans $\Xtilde$ et posons
$y=\somme_Qn_Q[\Orond_{C_Q}]$ dans $K_0(\Xtilde)^{(2)}$.

\th LEMME 20
\enonce 
On a $\;\psitilde(y)\circ\alpha=\vtilde\circ\ev\;$ en
tant qu'applications $\Atildeetoile\fleche\Z$.
\endth
\demonstration Cela résulte du lemme~19 --- qui en est le cas particulier où
$c$ est réduit à un point fermé --- par multiplicativité.\cqfd

Supposons de plus que $\deg c=0$.  Alors
$\ev(\Ktildeetoile)=\{1\}$ et $\alpha(\pi)=P$, donc $\ev$ induit
$\ev_0:\Atilde{}^\times\!/\,\Ktildeetoile\fleche \Ktildeetoile$
et $\alpha$ induit $\alpha_0:\Atilde{}^\times\!/\,\Ktildeetoile
\fleche T$, par passage aux quotients.  De plus, l'élément
$\psitilde(y)$ de $\Hom(\Z^S,\Z)$ appartient à $\Hom(T,\Z)$
(cf.~lemme~11). 

\th LEMME 21
\enonce 
On a $\;\psitilde(y)\circ\alpha_0=\vtilde\circ\ev_0:
\Atildeetoile/\Ktildeetoile\fleche\Z$. 
\endth
\demonstration C'est une conséquence immédiate du lemme~20.\cqfd

Soit $c$ un \ocycle\ de degré~0 sur $X_K$, notons $y\in K_0(X)^{(2)}$
la classe de son adhérence schématique dans $X$ et appliquons ce qui
précède à $\ctilde$, le \ocycle\ sur $X_\Ktilde$ obtenu à partir de
$c$ par changement de base.

\th LEMME 22
\enonce
On a
$\;\psi(y)\circ\alpha_0=\vtilde\circ\ev_0:
\Atildeetoile/\Ktildeetoile\fleche\Z$.  
\endth
\demonstration Il s'agit d'un cas particulier du lemme~21.\cqfd

Soient $U\subset X_K$ l'ouvert complémentaire du support $\Supp(c)$
de~$c$ et $\Utilde=U\times_K\Ktilde$.  Notons
$\cl_\Xtilde:\Zrond^1(\Utilde)\fleche\Pic\Xtilde$ l'application qui à
un diviseur de $\Utilde$ --- considéré comme un diviseur de
$\Xtilde$ --- associe sa classe dans $\Pic\Xtilde$.

\th LEMME 23
\enonce
Avec les notations précédentes, on a le carré cocartésien:
$$
\diagram{
\Atildeetoile/\Ktildeetoile&\droite{\div_\Utilde}&\Zrond^1(\Utilde)\cr
\vfl{-\alpha_0}{}{6mm}&&\vfl{}{\cl_\Xtilde}{6mm}\cr
T&\droite{}&\Pic\Xtilde.\cr
}
\advance\abovedisplayskip by-\baselineskip
\advance\belowdisplayskip by-2\baselineskip
$$
\endth
\demonstration Montrons que ce carré est commutatif. Pour tout
$g\in\Atildeetoile$, 
la classe du diviseur principal
$\div_\Xtilde(g)=\alpha(g)+\div_\Utilde(g)$ dans $\Pic\Xtilde$ est
nulle.  Comme cette classe vaut
$\alpha_0(g)+\cl_\Xtilde\div_\Utilde(g)$, on a
$-\alpha_0=\cl_\Xtilde\div_\Utilde$.  Le carré est cocartésien par la
définition de $\Pic\Xtilde$ comme quotient de la somme directe
$\Zrond^1(\Xtilde)=\Z^S+\Zrond^1(\Utilde)$.
\cqfd

On pose désormais $I=\Gal(\Kbar|\Ktilde)$.  Prenant les invariants
sous $I$ dans une suite exacte de $\Gal({\Kbar}|K)$-modules continus
discrets 
$$
\{1\}\fleche \Kbaretoile\fleche E\fleche \Pic \XKbar\fleche \zero,
\numeroter\newcount\extpicxkbarkbaretoile 
\global\extpicxkbarkbaretoile=\formuleno
$$
et en remarquant, d'une part, que le groupe $H^1(I,\Kbaretoile)$
s'annulle (Satz~90), et, d'autre part, que l'application canonique
$\Pic\XKtilde\fleche(\Pic \XKbar)^I$ est un isomorphisme 
\citer\ctsdescente({p.~386, (1.5.0)}) puisque $\Br\Ktilde=\zero$, 
on tire la suite exacte de $\Gal(\ktilde|k)$-modules
$$
\{1\}\fleche \Ktildeetoile\fleche E^I\fleche \Pic \XKtilde\fleche \zero.
\numeroter\newcount\invn \global\invn=\formuleno
$$
On a donc un homomorphisme qui envoie la classe de
\formule(\extpicxkbarkbaretoile) en celle de \formule(\invn):
$$
\formule(\extpicxkbarkbaretoile)\;\mapsto\;\formule(\invn) :
\Ext^1_K(\Pic \XKbar,\Kbaretoile)\droite{}
\Ext^1_k(\Pic \XKtilde,\Ktildeetoile).
\numeroter\newcount\homoextpicxkbarkprim
\global\homoextpicxkbarkprim=\formuleno
$$
Aussi, la valuation $\vtilde:\Ktildeetoile\fleche\Z$ prolongeant $v$
fournit un homomorphisme 
$
\Ext^1_k(\Pic \XKtilde,\Ktildeetoile)\fleche
\Ext^1_k(\Pic \XKtilde,\Z).
$
On note 
$$
\varepsilon:\Ext^1_K(\Pic\XKbar,\Kbaretoile)
\droite{} \Ext^1_k(\Pic\XKtilde,\Z) 
\numeroter\newcount\epsdef \global\epsdef=\formuleno
$$ 
l'application composée de \formule(\homoextpicxkbarkprim) suivie de
celle induite par $\vtilde$.
\goodbreak

Rappelons brièvement la définition de l'homomorphisme caractéristique
$$
\phi:\azeroo{X}\fleche \Ext^1_K(\Pic \XKbar,\Kbaretoile),
\numeroter \newcount\phidef \global\phidef=\formuleno
$$
de~\citer\ctssequel().  Soit $c$ un $0$-cycle de degré~$0$ sur $X_K$ ;
notons $[c]\in \azeroo{X_K}$ sa classe et posons $U=X_K-\Supp(c)$,
$\overline U=U\times_K \Kbar$.  Notons $\overline A$ l'anneau
semilocal de $\XKbar$ le long de $\cbar$, où $\cbar$ est le $0$-cycle
déduit de $c$ par changement de base.

Par définition, $\phi([c])\in\Ext^1_K(\Pic \XKbar,\Kbaretoile)$ est
la classe de l'extension 
de $\Pic\XKbar\,$ par $\Kbaretoile$ déduite de la suite exacte courte
$$
\{1\}\fleche \overline A{}^\times\!/\,\Kbaretoile\fleche 
\Zrond^1(\overline U) \fleche \Pic \XKbar \fleche\zero
\numeroter\newcount\extcan \global\extcan=\formuleno
$$
via l'homomorphisme
$\ev_0:\Abaretoile\!/\,\Kbaretoile\fleche\Kbaretoile$ d'évaluation en 
$\cbar$.

\th LEMME 24
\enonce
L'élément\/ $\varepsilon\phi([c])\in\Ext^1_k(\Pic\XKtilde,\Z)$ est la
classe de l'extension de\/ $\Pic\XKtilde$ par\/ $\Z$ déduite de la
suite exacte courte
$$
\{1\}\fleche \Atildeetoile\!/\,\Ktildeetoile\fleche 
\Zrond^1(\Utilde) \fleche \Pic \XKtilde \fleche\zero
\numeroter
$$
via\/ $\vtilde \circ \ev_0:\Atildeetoile \!/\,\Ktildeetoile\fleche\Z$ 
(\/l'évaluation en $\ctilde$ suivie de la valuation).
\endth
\newcount\extcantilde \global\extcantilde=\formuleno
\demonstration La longue suite exacte associée à la suite exacte courte
$$
\{1\}\fleche \Kbaretoile \fleche \Abaretoile \fleche \Abaretoile
/\Kbaretoile\fleche \{1\}
$$
de $I$-modules fournit un isomorphisme $\Atildeetoile/\Ktildeetoile
\fleche (\Abaretoile/\Kbaretoile)^I$, puisque $H^1(I,\Kbaretoile)=0$.
Il est alors clair que l'évaluation
$\ev_0:\Atildeetoile/\Ktildeetoile\fleche \Ktildeetoile$ en $\ctilde$
est la restrictions aux $I$-invariants de l'évaluation
$\ev_0:\Abaretoile/\Kbaretoile\fleche\Kbaretoile$ en $\cbar$.  Or la
suite~\formule(\extcantilde) s'obtient en prenant les $I$-invariants
dans la suite~\formule(\extcan) car l'application canonique
$\Pic\XKtilde\fleche(\Pic \XKbar)^I$ est un isomorphisme
\citer\ctsdescente({p.~386, (1.5.0)}) puisque $\Br\Ktilde=\zero$.\cqfd
\th LEMME 25
\enonce
L'élément\/ $\deltabar\gamma_0([c])\in\Ext^1_k(\Pic\XKtilde,\Z)$ est
représenté par l'extension déduite de la
suite exacte \formule(\extcantilde) via $\;\psi(y)\circ\alpha_0
:\Atildeetoile \!/\,\Ktildeetoile\fleche\Z$, où $y\in K_0(X)^{(2)}$
est un relèvement de $[c]\in\azeroo{X_K}$.
\endth
\demonstration  L'élément\/ $\deltabar\gamma_0([c])$ est la classe de
l'extension déduite de la suite \formule(\suitepic) via
$\;-\psi(y):T\fleche\Z$ (voir \citer\bourbaki(p.~125) pour le signe).
À son tour, \formule(\suitepic) est déduite de la suite
\formule(\extcantilde) via
$\;-\alpha_0:\Atildeetoile/\Ktildeetoile\fleche T$ (lemme~23).\cqfd

\th PROPOSITION 4
\enonce 
En tant qu'applications 
$\azeroo{X_K}\fleche \Ext^1_k(\Pic \XKtilde,\Z)$, on a\ \ 
$\varepsilon\circ\phi=\deltabar\circ\gamma_0$.
\endth
\demonstration  Soit $c$ un 0-cycle de degré~0 sur $X_K$.  Vus les
lemmes~24 et~25, il suffit de montrer qu'il existe un relèvement
$y\in K_0(X)^{(2)}$ de $[c]\in\azeroo{X_K}$ tel que
$\psi(y)\circ\alpha_0=\vtilde \circ \ev_0$ en tant qu'applications
$\Atildeetoile \!/\,\Ktildeetoile\fleche\Z$.  Or le lemme~22 affirme
précisément que c'est le cas lorsqu'on prend pour $y$ la classe de
l'adhérence schématique de $c$ dans $X$.\cqfd

\section{Conséquences et compléments}

{\it On suppose désormais que le groupe commutatif $\;\Pic\XKbar\;$
est libre de type fini} --- c'est le cas, par exemple, si $X_\Kbar\,$
est birationnel à $\P_2$ --- {\it et que l'application naturelle\/
$\Pic\XKtilde\fleche\Pic\XKbar\;$ est un isomorphisme.}

\th LEMME 26
\enonce 
L'application\/ $\varepsilon:\Ext^1_K(\Pic\XKbar,\Kbaretoile)
\fleche \Ext^1_k(\Pic\XKtilde,\Z)$ \formule(\epsdef) est alors un
isomorphisme.
\endth
\demonstration L'application $\varepsilon$ est composée de
\formule(\homoextpicxkbarkprim) et l'application
$$
\Ext^1_k(\Pic \XKtilde,\Ktildeetoile)\droite{}
\Ext^1_k(\Pic \XKtilde,\Z).
\numeroter\newcount\indvprim \global\indvprim=\formuleno
$$
induite par la valuation
$\vtilde:\Ktildeetoile\fleche\Z$. L'inclusion $\Ktilde\fleche\Kbar$
et l'identification $\Pic\XKtilde\fleche\Pic\XKbar\,$ 
fournissent une application
$$
\Ext^1_k(\Pic \XKtilde,\Ktildeetoile)\droite{}
\Ext^1_K(\Pic \XKbar,\Kbaretoile)
$$
réciproque de \formule(\homoextpicxkbarkprim).  Reste à faire voir que
\formule(\indvprim) est un isomorphisme.  L'uniformisante $\pi$ de $K$
fournit un scindage équivariant $\sigma:{\bf Z}\fleche\Ktildeetoile$
de la suite exacte
$$
\{1\}\fleche \ogothtildemult \fleche \Ktildeetoile
\hfl{\vtilde}{}{8mm} \Z \fleche \zero 
$$
de $\Gal(\ktilde|k)$-modules ; il induit par fonctorialité une
section de la suite 
$$
\zero\fleche\Ext_k^1(M,\ogothtildemult)\fleche\Ext_k^1(M,\Ktildeetoile)
\hfl{\formule(\indvprim)}{}{10mm}\Ext_k^1(M,\Z)\fleche \zero,
$$
avec $M=\Pic\XKtilde$, d'où l'exactitude de celle-ci.  Le groupe
$\Ext_k^1(M,\ogothtildemult)$ s'identifie à $H^1(\Spec\ogoth,D(M))$,
avec $D(M)=\Hom_k(M,\ogothtildemult)$, puisque $M$ est libre de type
fini
\citer\ctsdescente() ; à son tour, $H^1(\Spec\ogoth,D(M))$ s'identifie
à $H^1(k,D(M)_k)$
\citer\sgaiii(p.~401), qui est nul puisque le corps résiduel $k$ est
fini \citer\serrecg(p.~III.11), d'où le fait que \formule(\indvprim)
est un isomorphisme.
\cqfd

{\it On suppose de plus que
$\phi:\azeroo{X_K}\fleche\Ext^1_K(\Pic\XKbar,\Kbaretoile)$
\formule(\phidef) est injectif} --- hypothèse vérifiée lorsque 
$X_\Kbar\,$ est birationnelle à $\P_2$ \citer\hilbertxc().

\th THÉORÈME 1
\enonce
L'application 
$\gamma_0:\azeroo{X_K}\fleche B(X)_0$~\formule(\gammaodef) est alors un
isomorphisme.
\endth
\demonstration En effet, considérons le diagramme commutatif (prop.~4)
$$
\diagram{
\azeroo{X_K}&\hfl{\phi}{}{6mm}&\Ext^1_K(\Pic\XKbar,\Kbaretoile)\cr
\vfl{\gamma_0}{}{6mm}&&\vfl{}{\varepsilon}{6mm}\cr
B(X)_0&\hfl{\deltabar}{}{6mm}&\Ext^1_k(\Pic\XKtilde,\Z)\cr
}
\advance\abovedisplayskip by-\baselineskip
\advance\belowdisplayskip by-\baselineskip
$$
dans lequel, vues les hypothèses en vigueur, $\phi$ est injectif et
$\varepsilon$ est un isomorphisme (lemme~26).  Ceci montre que
$\gamma_0$ est injectif.  Or la prop.~3 affirme que $\gamma_0$ est
surjectif.
\cqfd

\th THÉORÈME 2
\enonce
Sous les même hypothèses, l'application de spécialisation
$\gamma:\azero{X_K}\fleche B(X)$~\formule(\gammadef) est un
isomorphisme. 
\endth
\demonstration En effet, considérons le diagramme comutatif
(lemme~15)
$$
\diagram{
\zero&\fleche&\azeroo{X_K}&\fleche&\azero{X_K}&\hfl{\deg}{}{6mm}&\Z&\cr 
&&\vfl{\gamma_0}{}{6mm}&&\vfl{\gamma}{}{6mm}&&\pafl{}{}{6mm}\cr
\zero&\fleche&B(X)_0&\fleche&B(X)&\hfl{\xibar}{}{6mm}&\Z&;\cr
}
\advance\abovedisplayskip by-\baselineskip
\advance\belowdisplayskip by-\baselineskip
$$
comme $\gamma_0$ est un isomorphisme (th.~1) et $\gamma$ est surjectif
(prop.~2), il résulte que $\gamma$ est bijectif.\cqfd

\th THÉORÈME 3 
\enonce
Avec les même hypothèses, si\/ $X_\ktilde$ est irréductible, alors\/
$\deg:\azero{X_K}\fleche\Z$ est un isomorphisme.
\endth
\demonstration Comme $\Card S=1$, l'application
$\xi:\Hom_k(\Z^S,\Z)\fleche\Z$ est un isomorphisme et
$\psi(i_*(K_0(X_k)^{(1)})=\zero$ (lemme~5).  Il en résulte que
$\xibar:B(X)\fleche\Z$ est un isomorphisme ; on utilise le lemme~15 et
le th.~2 pour conclure.\cqfd

\th THÉORÈME 4
\enonce
Supposons en outre que les composantes irréductibles (rendues
réduites) de $X_k$ sont lisses sur $k$.  L'application $\psi$
\formule(\psidef) induit alors un isomorphisme de $A_0(X_K)$ avec le
conoyau de l'homomorphisme
$$
\bigoplus_Y \Pic Y \fleche \Hom_k(\Z^S,\Z)
\quad(Y:\ \hbox{composantes irréductibles de\ } X_k)
\numeroter
$$
qui à une courbe $C$ sur un tel\/ $Y$ et une composante
irréductible (rendue réduite) $Z$ de $X_{\ktilde}$ associe le
nombre
$\deg_{\Ctilde}(\Orond_{\Xtilde}(Z)|_{\Ctilde})$~\formule(\lenombre),
où $\Ctilde=C\times_k \ktilde$. 
\endth
\newcount\caslisse \global\caslisse=\formuleno
\demonstration  Comme les $Y$ sont lisses sur $k$, on a  des
isomorphismes  
$$\Pic Y\droite{} {K^0(Y)^{(1)}\over K^0(Y)^{(2)}}\droite{}
{K_0(Y)^{(1)}\over K_0(Y)^{(2)}} 
$$
(cf.~\citer\manin(p.~45) et la démonstration du lemme~1).  On en
déduit que l'application
$\oplus_Y\Pic Y\fleche K_0(X_k)^{(1)}/K_0(X_k)^{(2)}$ est surjective ;
comme \formule(\caslisse) se factorise à travers celle-ci, son image
est égale à $\psi(i_*(K_0(X_k)^{(1)})$.  Le conoyau de
\formule(\caslisse) est donc égal à $B(X)$~\formule(\gammadef) ; le
th.~2 affirme que l'application $\gamma:\azero{X_K}\fleche B(X)$
induite par $\psi$ est un isomorphisme.\cqfd

Pour calculer ces nombres, il est souvent commode de se rappeler que
$\somme_Z m_Z\deg_{\Ctilde}(\Orond_{\Xtilde}(Z)|_{\Ctilde})=0$ (lemme~9).

{\it Exemple}\pointir Supposons que le nombre premier $p$ est impair.
Soient $d\in\ogothetoile$ non carré et
$\alpha,\beta,\gamma,\delta\in\Ketoile$ des éléments dont les
valuations satisfont $v(\alpha)=v(\beta)+1$ et
$v(\delta)=v(\gamma)+2$.  Dans l'espace $\P_{4,K}$ de
coordonnées $\def\sep{\,\colon}\,r\sep s\sep t\sep u\sep w$,
considérons la surface :
$$
\cases{\alpha(dr^2-s^2)=\beta(t-u)(t+w),\cr
\gamma(dr^2-t^2)=\,\delta(s+u)(s+w).\cr}
\numeroter\newcount\modelxo \global\modelxo=\formuleno 
$$
Soit $\pi$ une uniformisante de $K$ ; on
peut supposer que $\alpha=\pi$, $\delta=\pi^2$ et que
$\beta,\gamma\in\ogothetoile$.  Soit $X_0$ la $\ogoth$-surface
projective définie dans $\P_{4,\ogoth}$ par le système
\formule(\modelxo). La fibre fermée $X_{0,k}$ possède deux composantes
irréductibles, $A_0$ et $B_0$.  Éclatant la réunion des lieux
singuliers de $A_0$ et $B_0$, on obtient un $\ogoth$-schéma $X_1$ ; sa
fibre fermée $X_{1,k}$ possède quatre composantes irréductibles : les
transformées strictes $A_1$ et $B_1$ de $A_0$ et $B_0$, ainsi que deux
composantes exceptionnelles $C_1$ et $D_1$.  Les seuls singularités du
schéma $X_1$ sont, à part le point d'intersection $M_1$ de ces quatre
composantes, deux points $S_1$ et $R_1$. Éclatant ces trois points, on
obtient une $\ogoth$-surface $X$ qui est régulière, projective et
plate.  Sa fibre fermée $X_k$ possède sept composantes irréductibles :
les quatre transformées strictes $A$, $B$, $C$, $D$ et les trois
composantes exceptionnelles $R$, $S$, $M$.  Elles sont toutes de
multiplicité~1 sauf\/ $M$, qui est de multiplicité~2 ; rendues
réduites, elles sont lisses sur~$k$.  Seules $C$ et $D$ sont
absolument irréductibles ; les cinq autres possèdent deux composantes
irréductibles sur~$\ktilde$.

La fibre générique $X_K$ est $\Kprim$-birationnelle à $\P_{2,\Kprim}$,
où $\Kprim=K(\!\sqrt d)$ ; l'application $\Pic\XKtilde\fleche\XKbar\,$
est un isomorphisme.  Les groupes $\Pic Y$ des sept composantes
irréductibles $Y$ de $X_k$ sont faciles à calculer.  Appliquons le
th.~4 pour calculer $\azero{X_K}\,$ ; dans la base canonique de
$\Hom_k(\Z^S,\Z)$ (indexée par les composantes irréductibles $Y$ de
$X_k$) et pour un certain choix de générateurs des $\Pic Y$, l'image
de \formule(\caslisse) est engendré par les dix colonnes de la matrice 
{\catcode`\+=\active
\def+{\hbox{\phantom{$-$}}}
$$\def\\#1{{\goth #1}}
\bordermatrix{&\cr
A&-2&-1&-1&-2&+1&+1\cr
B&+1&  &  &  &  &-2&-1&-1&-2&+1\cr
C&  &  &  &  &  &+2&  &  &  &-2\cr
D&+2&  &  &  &-2\cr
R&  &+1&  &  &  &  &+1\cr
S&  &  &+1&  &  &  &  &+1\cr
M&  &  &  &+1&  &  &  &  &+1\cr
} 
$$}
\noindent
On en déduit que $\azeroo{X_K}$ est isomorphe à $\Z/2\Z$
(cf.~\citer\coomuzero(th.~4.5)) ; voir~\citer\delpezzo() pour plus
de détails.  

Il convient de remaquer que lorsque $d\in\Ketoile$ n'appartient plus à
$\ogothetoile$, l'application $\Pic\XKtilde\fleche\Pic\XKbar\;$ n'est
plus un isomorphisme pour la $K$-surface définie par
\formule(\modelxo).  Cependant, son groupe de Chow est
calculé dans \citer\chatelet() par une autre méthode~; celle-ci
n'exige plus l'existence d'une $\ogoth$-surface régulier, projectif
et plat de fibre générique \formule(\modelxo).

Rappelons pour terminer que lorsque $\XKbar\,$ est birationnel à $\P_2$,
le groupe commutatif $\Pic\XKbar\,$ est libre de type fini et que
l'homomorphisme caractéristique \formule(\phidef) est injectif
\citer\hilbertxc(), et, par suite, les théorèmes~1 à~4 s'appliquent
dès que $\Pic\XKtilde\fleche\Pic\XKbar\,$ est un isomorphisme.  Pour
d'autres exemples, voir \citer\barlowsimply(),
\citer\barlowrational().

\unvbox\bibbox
\bye